\documentclass[12pt]{article}
\usepackage{amssymb}
\title{Formules pour les nombres premiers}
\author{A.Balan \footnote{balan@paris7.jussieu.fr}}
\begin{document}
\maketitle
\abstract{The distribution of prime numbers is here considered. We show
 a formula for $li^{-1}$ and we study the
 $\pi(x)$ function and Riemann's hypothesis.}
\section{Introduction}
Il s'agit ici d'\'etudier la distribution des nombres premiers. En th\'eorie des nombres, on d\'efinit la fonction $\pi(x)$ comme \'etant
 le cardinal des nombres premiers plus petits que $x$ et on montre
 des estimations pour cette fonction, la premi\`ere \'etant $x/ln(x)$ \cite{H}.
 Une approximation de $\pi(x)$ avec le logarithme int\'egrale de Gauss 
$li(x)=\int_2^x \frac {1}{ln(t)} dt$ peut permettre de prouver 
l'hypoth\`ese de Riemann qui admet par ailleurs de multiples
 formes et qui se trouve impliquer de nombreux r\'esultats 
 en th\'eorie des nombres \cite{T}, sur les alg\`ebres d'op\'erateurs
 et autres domaines des math\'ematiques. La question est de localiser
 les z\'eros de la fonction zeta $\zeta(s)$ de Riemann sur la droite $Re(s)=1/2$. On sait que les z\'eros sont sym\'etriques par rapport
 \`a cette droite \`a cause de l'\'equation fonctionelle de $\zeta(s)$
 \cite{Ha}. Les approximations les meilleurs
 de la fonction $\pi(x)$ sont du genre $O(xexp(-cln(x)^{1/2}))$
 \cite{H}. 
 Il est possible de d\'emontrer qu'une
 approximation suffisante de type $O(x^{1/2+\epsilon})$
 pour $\pi(x)$ implique HR \cite{T}. Dans
 un premier temps, on calcule la fonction r\'eciproque de l'int\'egrale,
 ce qui permet d'approcher les nombres premiers par une formule. Ensuite,
 on red\'emontre qu'une bonne approximation de $\pi(x))$ implique HR.
 ¨Puis, par une autre m\'ethode, on montre qu'un prolongement
 analytique implique HR.
\section{La fonction r\'eciproque de $\int dx/ln(x)$}
On consid\`ere la fonction $\pi(x) \sim \int dx/ln(x)$. Il faut donc
 prendre en compte la fonction r\'eciproque
\begin{equation}
f(x)=(\int^x_2 \frac {1}{ln(t)}dt)^{-1}(x)
\end{equation}
Calculons maintenant la d\'eriv\'ee de $f$
\begin{equation}
f'(x)= \frac {1}{\frac {1}{ln(x)} \circ f}
\end{equation}
La d\'eriv\'ee est donc
\begin{equation}
f'= ln(f)
\end{equation}
On pose alors $f=e^g$ et on a 
\begin{equation}
g' e^g = g 
\end{equation}
On suppose $g^{(n)}=e^{-ng}P_n(g)$, on a alors une relation de r\'ecurrence sur les polynomes $P_n$ qui est
\begin{equation}
P_{n+1}(x)= x(-n P_{n}(x)+P'_n(x))
\end{equation}
\section{Une \'equation diff\'erentielle}
On pose maintenant
\begin{equation} 
l(x,y) = \sum_n P_n (x) {y^n/n!}
\end{equation} 
L'\'equation diff\'erentielle en $f$ est donc
\begin{equation}
l(x,y)= P_0(x) - xy l(x,y)+ x \partial_x \sum_n P_n(x) y^n/(n+1)!
\end{equation}
L'\'equation devient  alors
\begin{equation}
(1+xy)l(x,y)=P_0(x) + x \partial_x \int_0^y l(x,t)dt
\end{equation}
Cela donne pour $h= \int_0^y l(x,t)dt$,
\begin{equation}
(1+xy) \partial_y h= P_0 (x)+ x \partial_x h
\end{equation}
Et donc
\begin{equation}
(1+e^u y) \partial_y h = P_0(e ^u) +\partial_u h
\end{equation}
$u=ln(x),y=y$
$$
\partial_{y-u} h + e^u y \partial_y h = e^u
$$
$t=y-u, y=y$
$$
\partial_t h + e^{-t} y e^y \partial_y h = e^{y-t}
$$
$$
\partial_{e^{-t}} h + \partial_{\int 1/(ye^y)} h = 
e^t P_0(e^{y-t})
$$
$v=e^{-t}, w= -\int_1 1/(ye^y)$
$$
\partial_{v-w} h = e^{y(w)}
$$
$r=v-w, v=v$
On int\`egre alors et on trouve $h$ et donc une formule pour le $n$-i\`eme nombre premier en
 appliquant la formule de Taylor \`a $ln(f)$.
$$
\partial_r h = e^{y(v-r)}
$$
$$
h=\int^r_{\infty} e^{y(v-r)} dr
$$
Cela donne donc
$$
h(x,y)= \int^{e^{-y}x-\int_1 1/(ye^y)}_{\infty}(e^{y(-r+e^{-y}x)})dr
$$
Cela donne alors
$$
g(x+y)= \sum_n g^{(n)}(x)y^n/n! = \sum_n P_n(g(x)) (e^{-g(x)}y)^n/n! 
$$
$$
f(x+y)=e^{g(x+y)}= e^{l(g(x),e^{-g(x)}y)}=e^{\partial_y h(g(x),e^{-g(x)}y)}
$$
\section{L'hypoth\`ese de Riemann}
On fait ici quelques calculs avec $\zeta$. On montre qu'une 
 bonne approximation de $\pi(x)$ implique HR.
On d\'efinit la fonction $\zeta_p$ comme
\begin{equation}
\zeta_p (s) = \sum_{p,premier} \frac {1}{p^s}
\end{equation}
On a alors \cite{Bo}
\begin{equation}
\frac {\zeta'(s)}{\zeta(s)} = \sum_{p,premier} \frac {-ln(p)}{p^s (1-
 \frac {1}{p^s})}=\sum_{p,premier} \frac {-ln(p)}{p^s-1}=
\sum_{k \geq 1} \zeta_p'(ks)
\end{equation}
En effet, on a
$$
\frac{1}{1-1/p^s}= \sum_{k \geq 0} 1/p^{ks} 
$$
Et aussi, comme la convergence est en norme, on
peut inverser les signes somme.
$$
\sum_p -ln(p)/p^s \sum_k 1/p^{ks}= \sum_k \sum_p -ln(p)/p^{(k+1)s}
$$
La fonction $\zeta$ de Riemann ne peut avoir de z\'ero pour $Re(s) \succ \frac {1}{2}$
 si la fonction $\zeta_p'$ se prolonge analytiquement sur cette bande. En effet
\begin{equation}
\frac {\zeta'(s)}{\zeta(s)}-\zeta_p'(s)= \sum_{p,premier}
\frac {-ln(p)}{p^s(p^s-1)}
\end{equation}
La somme converge donc pour $Re(s) \succ \frac {1}{2}$.
Il faut alors comparer la fonction $\zeta_p$ avec la s\'erie suivante
\begin{equation}
\sum_n \frac {1}{((\int dx/ln(x))^{-1}(n))^s}
\end{equation}
Pour comparer les deux s\'eries, on forme la diff\'erence
\begin{equation}
\zeta_p (s) - \sum_n \frac {1}{((\int dx/ln(x))^{-1}(n))^s}
=
\end{equation}
$$
\sum_n [1-(1-(1-\frac {p_n}{(\int dx/ln(x))^{-1}(n)}))^s]/[p_n^s]
$$
On prend alors les normes pour montrer une convergence uniforme.
\begin{equation}
\sum_n |s|
|1- \frac {p_n}{(\int dx/ln(x))^{-1}(n)}| \frac {1}{p_n^{Re(s)}} +
\end{equation}
$$
\sum_n |s| |1- \frac {(\int dx/ln(x))^{-1}(n)}{p_n}| \frac {1}{(\int dx/ln(x))^{-1}(n) ^{Re(s)}}
$$
Il faut donc \'etudier les s\'eries suivantes
\begin{equation}
\sum_n \frac {1}{(\int dx/ln(x))^{-1}({n})^{Re(s)+1}}|(\int dx/ln(x))^{-1}({n})-p_n|
\end{equation}
$$
\sum_n \frac{1}{p_n^{Re(s)+1}}|(\int dx/ln(x))^{-1}({n})-p_n|
$$
Pour montrer que l'on peut prolonger analytiquement, il faut donc 
consid\'erer la diff\'erence 
\begin{equation}
|f(n)-f(\int^{p_n}_2 dx/ln(x))| \leq f'(c_n)|\pi(p_n)-\int^{p_n}_2 dx/ln(x)|
 \leq
\end{equation} 
$$
Cn^{1/2+\epsilon}
$$
Les s\'eries sont donc du type
\begin{equation}
\sum_n \frac{1}{(n)^{Re(s)+1/2-\epsilon}}
\end{equation}
Maintenant, on compare
\begin{equation}
\sum_n 1/f(n)^s , \int_0^{\infty} dt/ f(t)^s
\end{equation}
On fait un changment de variable dans l'int\'egrale $u=f(t)$ et on d\'erive par rapport \`a $s$ pour obtenir
\begin{equation}
\int_2^{\infty} (1/ln(u))(1/u^s)du =- \int_{+\infty}^s (2^{1-r}/1-r)dr 
\end{equation}
Elle est donc prolongeable analytiquement. La diff\'erence est de plus
\begin{equation}
\int_0^{\infty} |[1/f(t)^s-1/f(E(t))^s]| dt \leq \sum_n \int_n^{n+1}
\frac {|f(E(t))^s-f(t)^s]|}{[f(t)f(E(t))]^{Re(s)}} dt \leq
\end{equation}
$$
\sum_n |s| \int_n^{n+1}
 f'(c_t)(t-E(t))/[f(E(t))]^{Re(s)+1} dt \leq
$$
$$
C 1/(n)^{Re(s)+1/2-\epsilon}
$$
La diff\'erence converge en norme, ce qui va entra\^\i ner que l'on peut prolonger analytiquement sous cette condition.
\section{Les z\'eros de la fonction $\zeta$}
On a \cite{T}
\begin{equation}
\tilde \zeta(s)= \prod_p \frac {1}{1+ \frac {1}{p^s}}
\end{equation}
$\tilde \zeta$ est d\'efinie pour $Re(s) \succ 1$. Alors
\begin{equation}
\tilde \zeta (s) \zeta (s)= \zeta (2s)
\end{equation}
 Cette formule est v\'erif\'ee pour $Re(s) \succ 1/2$.
Cela entraine que si $\zeta(s)$ est nul, il faut que 
 soit $\zeta(2s)$ est nul, soit $\tilde \zeta(s)$ infini. Or pour $Re (s) \succ 1$, $\zeta(s)$ ne peut avoir de z\'ero.  On
 ne peut avoir de z\'ero si on peut prolonger analytiquement pour $\tilde \zeta(s)$. 

\bigskip

Si on avait un z\'ero $s_0$ pour zeta entre $1/2$ et $1$,
 on aurait un p\^ole pour $\zeta_p'$ et alors on se place en $s_0/2$ et on a
\begin{equation}
\zeta'(s_0/2) / \zeta (s_0/2)= \zeta_p'(s_0/2)+ \zeta_p'(s_0)+ \ldots
\end{equation}
\begin{equation}
\tilde \zeta'(s_0/2) / \tilde \zeta (s_0/2)= -\zeta_p'(s_0/2)+ \zeta_p'(s_0)+ \ldots
\end{equation}
Les autres termes sont convergents en dehors de $3/2s_0$. 
La premi\`ere \'egalit\'e donne qu'il y a un p\^ole pour
 $\zeta_p'(s_0/2)$ car on suppose par r\'ecurence que les
 z\'eros sont sur la droite $1/2$ pour $Im(s)< Im(s_0)$.
 On consid\`ere alors la seconde \'egalit\'e et donc un z\'ero double serait pour $\tilde \zeta$ en $s_0/2$~; ce qui contredit
 $\tilde \zeta(s) = \zeta(2s)/\zeta(s), \tilde \zeta (s_0/2)=\zeta(s_0)/\zeta(s_0/2)=0$ et d\'emontre HR.
\footnote{Il faut aussi prendre en compte les z\'eros 
multiples et les z\'eros pour $3/2 s_0$.}
\section{Derni\`eres formules}
 On a de plus
\begin{equation}
\frac {\tilde \zeta'(s)}{\tilde \zeta(s)}=
\sum_{k \geq 1} (-1)^k \zeta_p' (ks)
\end{equation}
On a aussi 
\begin{equation}
\zeta'(s)/\zeta(s) - \zeta'(2s)/\zeta(2s)= \sum_{k=2n+1} \zeta'_p(ks)
\end{equation}
Et aussi
\begin{equation}
\zeta'_p(s)= \zeta'(s)/\zeta(s)-\zeta'(2s)/\zeta(2s)-\zeta'(3s)/\zeta(3s)
+ \ldots
\end{equation}
\begin{equation}
\zeta_p(s)=\sum_n (\mu(n)/n) ln(\zeta(ns))
\end{equation}
De plus, une g\'en\'eralisation de l'\'equation fonctionelle de zeta est
\begin{equation}
(\sum_n 1/(n^2+an)^s) \Gamma_{\pi}(s)= (\sum_k cos(\pi ka)/k^{1-2s})\Gamma_{\pi}(1/2-s)
\end{equation}
\section{Remerciements}
Je remercie chaleureusement J.Kouneiher, C.Peschard et J.-L.Tu 
pour leurs conseils.

\end{document}